\documentclass[a4paper,12pt]{article}
 \usepackage{amssymb,latexsym,amsmath,amsthm}
\theoremstyle{definition}
\newtheorem{opr}{Definition}
 \theoremstyle{plain}
\newtheorem{thm}{Theorem}
\newtheorem{lem}{Lemma}
\newtheorem{utv}{Proposition}
\newtheorem*{sled}{Corollary}

\theoremstyle{remark}
 \newtheorem*{note}{Remark}
\newtheorem*{exam}{Example}

\DeclareMathOperator{\rank}{rk}

\newcommand{\bl}{\square}

\newcommand{\ma}{\mathcal{A}}
\newcommand{\ms}{\mathcal{S}}

\newcommand{\R}{\mathbb{R}}
\newcommand{\Z}{\mathbb{Z}}
\newcommand{\Q}{\mathbb{Q}}

\newcommand{\doc}{\emph{Proof. }}
\renewcommand{\leq}{\leqslant}
\renewcommand{\geq}{\geqslant}

\newcommand{\ph}{\varphi}
\newcommand{\itd}[2]{#1_1,\dotsc, #1_{#2}}

\newcommand{\Pm}{\mathbb{P}}

\newcommand{\La}{\Lambda}

\newcommand{\odo}[1]{=1,\dots, #1}

\newcommand{\Ph}{\varPhi}
\newcommand{\nil}{\mathfrak {N}_2(n,t)}
\newcommand{\nilk}{\mathfrak {N}_2^K(n,t)}
\newcommand{\nilo}{\mathfrak {N}_2(n_0,t_0)}
\newcommand{\nilg}{\mathcal{N}_2(n,t)}
\newcommand{\nilog}{\mathcal {N}_2(n_0,t_0)}

\numberwithin{equation}{section}

\begin{document}
\noindent\textsf{ MSC: 17B30, 20F18, 20F40, 14A10, 14M15\\ }

\begin{center}{\LARGE  \textsf{Generic properties of 2-step nilpotent Lie algebras and torsion-free groups}\\ }  
\end{center}

\begin{center}
 \textsf{\large Maria V. Milentyeva} 
\end{center}

{\footnotesize
  To define the notion of a generic property of finite dimensional 2-step nilpotent  Lie algebras  we use standard correspondence between such Lie algebras and points of an appropriate algebraic variety, where a negligible set  is one  contained in a proper Zariski-closed subset. We compute the  maximal dimension of an abelian subalgebra of a generic Lie algebra and  give a sufficient condition for a generic Lie algebra  to admit no surjective homomorphism onto a  non-abelian Lie algebra of a given dimension.
 Also we consider analogous questions for finitely generated torsion free  nilpotent groups of class 2.\\}

\section{Introduction and Results}

\subsection{Nilpotent Lie algebras of class 2}
Consider a finite-dimensional 2-step nilpotent Lie algebra $L$. Denote by $S$ the commutator subalgebra of $L$. Let $\itd{z}{t}$ be a basis of $S$, and let $V\subset L$ be such  a subspace that  $L= V \oplus S$. Then the product of two elements $x=\Bar{x}+\Bar{\Bar{x}}$ and  $y=\Bar{y}+\Bar{\Bar{y}}$ of $L$ with  $\Bar{x},\Bar{y}\in V$ and $ \Bar{\Bar{x}},\Bar{\Bar{y}}\in S$ has the form 
\begin{equation} \label{E:g2}
 [x,y]=\ph_1(\Bar{x},\Bar{y})z_1+ \dots+ \ph_t(\Bar{x},\Bar{y})z_t    
\end{equation}
for some $t$-tuple of alternating bilinear forms $\varPhi=\varPhi({L})=(\itd{\ph}{t})$ on $V$.
 
On the other hand, given vector spaces $S$ and $V$, a basis $\itd{z}t$ of $S$ and a $t$-tuple $\varPhi=(\itd \ph t)$ of alternating bilinear forms on $V$,  one can define the product of two elements of  ${L}={L}(\varPhi)=V \oplus S$ by~(\ref{E:g2}). Obviously, $L(\varPhi)$ is a 2-step nilpotent Lie algebra and $S$ is a central subalgebra of $L(\varPhi)$.

Let $\itd {z'} t$  be a new basis of $S$ and let $C=(c_{ij})$ be the transformation matrix   from old to new coordinates. Then~(\ref{E:g2}) can be written in the following form:
\begin{equation} \label{ea}
 [x,y]=\sum_{i=1}^t \ph_i(\Bar{x},\Bar{y})z_i= 
\sum_{i=1}^t \ph_i(\Bar{x},\Bar{y})\sum_{j=1}^t c_{ji} z'_j= 
 \sum_{j=1}^t( \sum_{i=1}^t \ph_i(\Bar{x},\Bar{y})c_{ji}) z'_j=
\sum_{i=1}^t \ph'_i(\Bar{x},\Bar{y}) z'_i,
\end{equation}
where 
\begin{equation} \label{eb}
\ph'_i =  \sum_{j=1}^t c_{ij}\ph_j,
\end{equation}
That is  $C^\bot$ is the transformation matrix to go from $\varPhi$ to  $\varPhi'=(\itd {\ph'} t)$.
This  gives the following properties of a $t$-tuple $\varPhi$.
 
\begin{utv} \label{ta} 
The alternating bilinear forms of  the $t$-tuple $\varPhi(L)$ corresponding to a Lie algebra $L$ are linearly independent. And conversely, if the forms of a $t$-tuple $\Ph$ are linearly independent and $L(\Ph)=V\oplus S$ is a corresponding Lie algebra, then $S=[L(\varPhi), L(\varPhi)]$.
\end{utv}
\doc 
Suppose, on the contrary, that for some Lie algebra $L$ the forms of $\Ph(L)$ are linearly dependent. Then there exists  non-trivial linear combination  $\alpha_1 \ph_1+ \dots \alpha_t \ph_t=0$. It follows from~(\ref{eb}) that, putting $c_{tj}=\alpha_j$, one can choose a basis $\itd {z'} t$ of $S$ in such a way that $\ph'_t\equiv 0$. This contradicts with the fact that $z'_t$ belongs to the commutator subalgebra $S$.

Conversely, if the commutator subalgebra  $[L(\varPhi), L(\varPhi)]$ is strictly less than $S$, then there exists a basis $\itd {z'} t$ of $S$ such that $\itd z k$, where $k< t$, generate $[L(\varPhi), L(\varPhi)]$. Then the corresponding bilinear form $\ph'_t$, which is a non-trivial linear combination of the forms of $\Ph$, identically equals zero.
~$\bl$

\begin{utv} \label{tb}
If the linear spans of $t$-tuples $\Ph$ and $\Ph'$ coincide, then the corresponding Lie algebras $L(\varPhi) $ and $ L(\varPhi')$ are isomorphic.
\end{utv}
 From now on we will assume that a ground field is infinite.
 It is appropriate to introduce some notation that will be used throughout the paper. We denote by
 $\nilk$ (or simply by $\nil$, if it is not misleading)  the set of all 2-step nilpotent Lie algebras $L$ over a field $K$ such that   $\dim L/[L,L]=n$ and $\dim [L,L]= t$.   Let  $B_n(K)$ be the set of all alternating bilinear forms on an $n$-dimensional space $V$ over $K$ and   ${\Pm}(B_n^t)$ be the projective space corresponding to the $t$-th direct power of $B_n(K)$. Put
$ M_{n,t}=\{(\ph_1: \dots :\ph_t)\in {\Pm}(B_n^t) \mid \itd \ph t
\mbox{ are lineary independent}\}.
$
$M_{n,t}$ is a Zariski-open subset of $ {\Pm}(B_n^t)$.

\begin{opr} \label{ty}
 We  say that  \emph{$\ma$ is a generic property of $\nil$  (a
generic Lie algebra of $ \nil$ has a property $\ma$)} if  the set of
points of $M_{n,t}$ corresponding to the algebras without the
property $\ma$ is contained in  some proper Zariski-closed subset. 
\end{opr}

\begin{note}
 Since $M_{n,t}$ is irreducible,  the dimension of any its
closed subset is strictly less than the dimension  of $M_{n,t}$. (see.~\cite{S88})
\end{note}

The following theorem, which holds true for any infinite field, gives a simple example of a generic property.
  
\begin{thm} \label{ti}
  For a
generic Lie algebra $ L \in \nil$ we have $\dim Z(L)=2$ if  $t=1$ and $n$ is odd, and $Z(L)=[L,L]$ otherwise.
 \end{thm}

Here $Z(L)$ is the center of a Lie algebra $L$.

The next theorem allows to compute the maximal dimension of an abelian subalgebra of a generic Lie algebra.

\begin{thm} \label{ts}
If a ground field $K$ is algebraically closed and $t>1$, then any
Lie algebra of $ \nilk$  contains a commutative subalgebra of
dimension $s=[\frac{2n+t^2+3t}{t+2} ]$. For any infinite ground field $K$
 a  generic Lie algebra of $\nilk$   doesn't have any
commutative subalgebra of dimension $s+1$.
\end{thm}

This theorem immediately follows from~\cite{M7}. However for the convenience of the reader we will give some details in Section~\ref{s6}. Also in this section we adduce an example of a Lie algebra from  $\mathfrak{N}^{\R}_2(4,3) $ without commutative subalgebras of dimension~5. Which shows that the condition for the ground field being algebraically closed can not be ignored. 

\begin{note}
The structure of Heisenberg Lie algebras, that is  2-step nilpotent Lie algebras with one-dimensional center, is well known. The dimension $m$ of such an algebra is always odd, and the dimension of a maximal commutative subalgebra is~$\frac{m+1}2$.
\end{note}

The main result of the present paper it the following theorem, which  holds true for an arbitrary infinite ground field.

\begin{opr} \label{tp}
 We say that a  nilpotent Lie algebra $L$ of class 2 \emph{ has a
property}
 $\ms$$(n_0,t_0)$ $(1\leq t_0\leq \frac{n_0(n_0-1)}2)$, if $L$ admits a surjective homomorphism onto a Lie algebra of $
 \nilo$.
\end{opr}

\begin{thm} \label{tf}
If the positive integers   $n,\ n_0,\ t$ and $t_0$ satisfy
the inequality
\begin{equation} \label{er}
t< \frac{n(n-1)}2-  \frac{nn_0}{t_0}  +\frac{n^2_0}{t_0}+t_0-
\frac{n_0(n_0-1)}2,
\end{equation}
then a generic Lie algebra of $\nil$ does not have the property
$\ms$$(n_0,t_0)$. If  $n,\ n_0,\ t$ and $ t_0$, where  $n\geq n_0 $,
satisfy the inequality
\begin{equation} \label{eq}
 t \geq \frac{n(n-1)}2-\frac{n_0(n_0-1)}2 +t_0,
\end{equation}
then  the property $\ms$$(n_0,t_0)$ is true on $\nil$.
\end{thm}

The following statement is an immediate corollary of Theorem~\ref{tf}.

\begin{sled} \label{tk}
A generic Lie algebra of $\nil$ does not admit a
surjective homomorphism onto a non-commutative Lie algebra of
dimension~$N< n$, if 
\begin{equation} \label{en}
t< \frac{n^2}2-n \left(N - \frac 1 2  \right) + \frac {N(N-1)} 2 +1.
\end{equation}
\end{sled}

\begin{note}
Obviously, if $n< N\leq n+t$, then any Lie algebra $L\in \nil$ admits a surjective homomorphism onto a non-cummutative Lie algebra of dimension~$N$. It is enaugh to take a 
quotient group of $L$ by any central subalgebra of dimension $n+t-N<t$.
\end{note}

Actually, if a ground field is algebraically closed, then the set of point of $M_{n,t}$ corresponding to  the algebras with property $\ms$$(n_0,t_0)$ forms a  closed subset  for any values of parameters (see Lemmas \ref{tc} and \ref{UT}). It means that there are only two possibilities: either  this subset  coincides with $M_{n,t}$ and all the algebras have the property $\ms$$(n_0,t_0)$, or this subset is proper and   almost all the algebras do not have this property.
Moreover, if all the allgebras of $\nil$ have the property $\ms$$(n_0,t_0)$, then it is also true for any Lie algebra $L \in \mathfrak N$$_2(n, t+k)$, where $k>0$, since for any $k$-dimensional subalgebra  $H < [L,L]$ we have $L/ H \in \nil$. Thus, we get the following.

\begin{utv} \label{tz}
 If a ground field is algebraically closed, then for any
positive integers $n,\ n_0$ and $t_0$, where $n\geq n_0$ and
$t_0\leq \frac{n_0(n_0-1)}2$, there exists an integer $C(n,n_0,t_0)$
such that for any integer $t$ satisfying the inequality
$$1\leq t < C(n,n_0,t_0) ,$$
a generic Lie algebra of $\nil$ doesn't have the property
$\ms$$(n_0,t_0)$, and  the property $\ms$$(n_0,t_0)$ is true on
$\nil$ if
$$C(n,n_0,t_0)\leq t \leq \frac{n(n-1)}2.$$
\end{utv}

The relations~(\ref{er}) and~(\ref{eq}) give upper and lower bounds for $C(n,n_0,t_0)$.
Computing  $C(n,n_0,t_0)$ precisely is  subject to further investigation.

\subsection{Finitely generated torsion-free nilpotent groups of class 2}

Now let us introduce the notion of a generic property for nilpotent groups in an analogous way.
Let $G$ be a  finitely generated 2-step nilpotent group without torsion. Denote by  $S=I(G')=\{x\in G\mid x^k\in G'\mbox { for some }k\}$ 
 the isolator of the commutator subgroup of $G$. Since the center of a torsion-free nilpotent group 
coincides with its isolator (see, for example,~\cite[Section 8]{B85}), $I(G')$ also lies in the center. Hence $S$ and the quotient group $G/S$ are free abelian.
Let $\itd b t$ and $\itd a n$ be  bases of $S$ and $G$ modulo $S$ respectively.
Then the elements of $G$ have the form
\begin{equation} \label{E:tt}
  a_1^{k_1}\dotsm a_n^{k_n}b_1^{l_1}\dotsm b_t^{l_t} , \mbox{ with }  k_i,l_j \in { \mathbb{Z}},\mbox{ for } \ i\odo n,\  j\odo t.
\end{equation}

Since the commutator map on a nilpotent group of class 2 is bilinear, we get
\[
   [x,y]=[a_1^{k_1}\dotsm a_n^{k_n}b_1^{l_1}\dotsm b_t^{l_t}, a_1^{p_1}\dotsm a_n^{p_n}b_1^{q_1}\dotsm b_t^{q_t}]=  
\] 
\begin{equation}  \label{E:i1}
=\prod_{i,j=1}^{n} [a_i,a_j]^{k_ip_j} = \prod_{l=1}^{t} b_l^{\sum_{i,j=1}^n  k_ip_j \ph_l(a_i,a_j)} = \prod_{l=1}^t b_l^{\ph_l(xS,yS)}.
\end{equation}
for some $t$-tuple of integer skew-symmetric bilinear forms $\varPhi(G)=\{\itd \ph t\}$ on $G/S$.

Conversely, let  $\varPhi=\{ \itd{\varphi}{t} \}$ be a $t$-tuple of integer skew-symmetric bilinear forms on a free abelian group $\mathbb{Z}^n$ with a basis $\itd{a}{n}$.
Consider the set $G(\varPhi)$ of formal products of the form~(\ref{E:tt}). Define the product of two elements of $G(\varPhi)$ by
\begin{align} \label{E:h1}
  a_1^{k_1}\dotsm a_n^{k_n}b_1^{l_1}\dotsm b_t^{l_t} \centerdot a_1^{p_1}\dotsm a_n^{p_n}b_1^{q_1}\dotsm b_t^{q_t}= \\
=a_1^{k_1+p_1}\dotsm a_n^{k_n+p_n}\prod_{k=1}^t b_k^{l_k+q_k+\sum_{i>j} k_i p_j \ph_k (a_i,a_j)}.   \notag
\end{align}

The product is well defined and it can be verified easily that the group axioms are satisfied. From~(\ref{E:h1})  it follows that $G(\varPhi)$ is a nilpotent torsion-free group of class 2, that the subgroup $S$ generated by $\itd{b}{t}$  is central and freely generated by $\itd{b}{t}$, and the factor group $G(\varPhi)/ S$ is abelian and freely generated by $a_1S,\dotsc, a_n S$. So we can naturally identify  it with $\Z^n$.
Using~(\ref{E:h1}), we get  
\begin{equation*}   
 [a_i,a_j]=\prod_{l=1}^{t} b_l^{ \ph_l(a_i,a_j) } \quad \mbox{ for  }\, i,j\odo n.
\end{equation*}
Further, since the commutator map is bilinear, we get that~(\ref{E:i1}) holds.

It is easy to see that  the   $t$-tuple of bilinear forms corresponding to a new basis of $S$ is defined by~(\ref{eb}),
 where the transformation matrix $C$ can be an arbitrary integer matrix with determinant 1.
 Therefore, as in the case of Lie algebras, we get the following property of $\Ph(G)$.
\\ \\
\textbf{Proposition~\ref{ta}$'$.} \emph{The alternating bilinear forms of  the $t$-tuple
$\varPhi(G)$ corresponding to a group $G$ are linearly independent.
And conversely, if the forms of  a $t$-tuple  $\varPhi$ are linearly
independent and $G(\Ph)$ is a corresponding 2-step nilpotent group,
then $S=I((G(\Ph))')$.
}\\

The proof is carried over   from Proposition~\ref{ta}. 

Let  $\mathcal{N}_2(n,t)$ be  the set of all 2-step nilpotent torsion free  groups $G$ with 
$\rank G'=t$ and $\rank G/ I(G')=n$.
Denote by $B_n(\mathbb{Z})$ the set of all integer bilinear forms on $\Z^n$. 
Put
$ P_{n,t}=\{(\ph_1, \dots ,\ph_t)\in B_n^t(Z) \mid \itd \ph t
\mbox{ are lineary independent}\}.
$
Given a basis $\itd a n$ of a vector space $V$ over $\Q$, 
we can identify the elements of $B_n(\Z)$ with elements of $B_n(\Q)$ represented in this basis by integer matrices.\\ \\
\textbf{Definition~\ref{ty}$'$.}
We say that \emph{a property $\ma$ of $\nilg$ is true
generically} if the set of points of $P_{n,t}$ corresponding to the
groups without property $\ma$ is contained in some proper Zariski-closed
subset of $B_n^t(\Q)$.\\

A nilpotent group $G\in \nilg$  and the Lie algebra $L=L(\Ph(G))\in \mathfrak{N}_2^{\Q}(n,t)$, defined by the same tuple of the forms, have quite similar properties. In particular, $G$ contains an abelian subgroup of rank $s$ if and only if $L$ has an $s$-dimensional commutative subalgebra. $G$ admits a surjective homomorphism onto a group of $\nilog$ if and only if $L$ has the property $\ms$$(n_0, t_0)$. Also we have  $\rank Z(G)= \dim Z(L)$ and $\rank G'=\dim [L,L]$. That is why the analogues of Theorems~\ref{ti}-\ref{tf} holds true for the nilpotent groups. They do not need to be proved separately. 
  This can be explained by the following arguments.

Denote by $\sqrt{G}$ the Malcev completion of $G\in \nilg$, that is the smallest complete nilpotent torsion free group containing~$G$ (see, for example,~\cite[Section 8]{B85}). Notice that   $\sqrt G$ can be considered as the group of elements of the form~(\ref{E:tt}) with multiplication~(\ref{E:h1}), where $(\itd \ph t)= \Phi(G)$ and $k_i,$$ l_j,$$p_r,$$ q_s \in \Q$.

On the other hand, it is well known (see~\cite{B85}) that in every nilpotent Lie $\Q$-algebra $L$ one can define multiplication "$\circ$"  by Campbell-Hausdorff formula  in such  a way that $L$ is a torsion-free complete group, say $L^\circ$, of the same nilpotency class   with respect to "$\circ$". (In the case of a 2-step nilpotent Lie algebra Campbell-Hausdorff formula has the form    $x \circ y = x+y+ \frac 1 2 [x,y]$.) The functor  $\Gamma: L \to L^\circ, (f:L_1\to L_2)\to(f:L^\circ_1\to L^\circ_2) $
provides an isomorphism from category of the nilpotent class $s$ Lie $\Q$ algebras to the category of the nilpotent class $s$ torsion-free complete groups. It is easy to see that 
$\sqrt G\simeq L(\Phi(G))^\circ$.
\\ \\
\textbf{Theorem~\ref{ti}$'$.} \emph{For a generic group $G\in \nilg$ we have $\rank Z(G)=2$ if  $t=1$ and  $n$ is odd, and $Z(G)=I(G')$ otherwise.}\\  \\
\textbf{Theorem~\ref{ts}$'$.} \emph{A generic group of $\nilg$ doesn't contain any abelian
subgroup of rank $\left[\frac{2n+t^2+3t}{t+2}\right]+1$.}\\

Also this theorem follows from~\cite{M4}.\\ \\ 
 \textbf{Theorem~\ref{tf}$'$.} \emph{Suppose that the positive integers $n,\ n_0, \ t,$ and
$t_0$ satisfy the inequality
\begin{equation*}  
t< \frac{n(n-1)}2-  \frac{nn_0}{t_0} +\frac{n^2_0}{t_0}+t_0-
\frac{n_0(n_0-1)}2.
\end{equation*}
Then a generic group of $\nilg$ does not  admit a surjective
homomorphism onto a group of $\nilog $.}

\emph{If  $n,\ n_0,\ t$ and $ t_0$, where  $n\geq n_0 $, satisfy the
inequality
\begin{equation*}  
 t \geq \frac{n(n-1)}2-\frac{n_0(n_0-1)}2 +t_0,
\end{equation*}
then every group of $\nilg$ admits a surjective homomorphism onto a
group of $\nilog$. }

\begin{sled}
A generic group of $\nilg$ does not admit a surjective
homomorphism on a non-abelian group of polycyclic rank~$N
< n$, if
\begin{equation*} 
t< \frac{n^2}2-n \left(N - \frac 1 2  \right) + \frac {N(N-1)} 2 +1  .
\end{equation*}
\end{sled}

\section{ Preliminaries on the algebraic varieties } \label{S11}

We  recall that $X \subset {\Pm}^n$ is a \emph{closed in Zariski topology subset} if it consists of all points at which a finite number of homogeneous  polynomials with coefficients in $K$ vanishes.
This topology induces Zariski topology on any subset of  the projective space ${\Pm}^n$. A closed subset of ${\Pm}^n$ is    a \emph{ projective variety}, and an open subset of a projective variety is   a \emph{quasiprojective variety}. A nonempty   set $X$ is called \emph{irreducible} if it cannot be written as the union of two proper closed subsets.

Further, let   $f:X\to\Pm$$^m$ be a map of a  quasiprojective variety $X\subset\Pm$$^n $ to a projective space $\Pm$$^m$. This map is \emph{regular} if for every point $x_0 \in X$ there exists a neighbourhood $U \ni x_0$ such that the map $f: U \to \Pm$$^m$ is given by an $(m+1)$-tuple $(F_0:\ldots:F_m)$ of homogeneous polynomials of the same degree in the homogeneous coordinates of $x \in {\Pm}^n$, and $F_i(x_0)\neq 0$ for at least one $i$.

We will  use the following properties of the \emph{dimension} of a quasiprojective variety.
\begin{itemize}
\item [-] The dimension of $\Pm$$^n$ is equal to $n$.
\item[-] 
 If $X$ is an irreducible variety and $U\subset X$ is open, then 
  $\dim U=\dim X$. 
\item [-]
 The dimension of a reducible variety is the maximum of the dimension of its irreducible components.
\end{itemize}
 
For more details see ~\cite{S88}. We will need the following propositions which are also to be found in~\cite{S88}.

\begin{utv} \label{T:a2}
If  ${X=\cup U_\alpha}$ with open sets $U_\alpha$, and   $Y\cap U_\alpha$ is closed in  $U_\alpha$ for each $U_\alpha$, then $Y$ is closed in $X$.
\end{utv}

\begin{utv} \label{tl} 
(\normalfont{see} \cite[Theorem~11.12]{H05})
Let $f: X \to \Pm$$^n$ be a regular map of a quasiprojective variety $X$, and let $Y$ be the closure of $f(X)$. For each point $p\in X$, we denote by $X_p= f^{-1}(f(p)) \subset X$ the fiber of $f$ containing $p$, and by $\dim_p X_p$ denote the local dimension of $X_p$ at the point $p$, that is, the maximal  dimensions of an irreducible component of $X_p$, containing $p$.
Then the set of all points $p \in X$ such that $\dim_p X_p \geq m$ is closed in $X$ for any $m$.
More over, if $X_0$ is an irreducible component of $X$ and $Y_0 \subset Y$ is the closure of the image $f(X_0)$, then
\begin{equation} \label{fez}
\dim X_0= \dim Y_0 + \min_{p\in X_0} \dim_p X_p. 
\end{equation}
\end{utv}

\begin{utv} \label{T:A2}
Let  $f:X\to Y$ be a regular map between projective varieties, with ${f(X)=Y}$. Suppose that $Y$ is irreducible, and that all the fibers  $f^{-1}(y)$ for $y \in Y$  are irreducible and of the same dimension. Then $X$ is irreducible.
\end{utv}

The following  examples of projective varieties are important for our goals: \emph{a direct product} of projective spaces, a \emph{Grassmannian} and a \emph{Shubert cell}.

Let ${\mathbb P}^n ,\ {\mathbb P}^m$ be projective spaces having homogeneous coordinates ${(u_0: \dotsc: u_n)}$  and 
 ${(v_0: \dotsc: v_m)} $ respectively. Then the set   ${\mathbb P}^n \times {\mathbb P}^m$ of pairs $(x,y)$ with $x\in {\Pm}^n$ and $y\in {\Pm}^m$ is naturally embedded as a closed set into the projective space  ${\mathbb P}^{(n+1)(m+1) -1}$ with homogeneous coordinates $w_{ij} $ by the rule $w_{ij}(u_0: \dotsc: u_n;v_0: \dotsc: v_m)=u_iv_j$.  Thus there is a topology on ${\mathbb P}^n \times {\mathbb P}^m$, induced by the Zariski topology on $\Pm^{(n+1)(m+1) -1}$.
 
\begin{utv} \label{T:e2}
  A subset $X \subset {\mathbb P}^n \times {\mathbb P}^m\subset \Pm$$^N$ is a closed algebraic subvariety if and only if it is given by a system of equations  
\[
 G_i(u_0: \dotsc: u_n;v_0: \dotsc: v_m)=0 \qquad (i \odo t)
\]
homogeneous in each set of variables $u_i$ and $v_j$.
\end{utv}

\begin{utv} \label{T:u2}
   Consider an  $n$-dimensional vector space  $V$ with a basis  $\{e_1, \dotsc ,e_n\}$. Let   $U$ be a  $k$-dimensional subspace of  $V$  with a basis $\{f_1, \dotsc ,f_k\}$. To $U$  we  assign the point  $P(U)$ of the projective space $\mathbb{P}$$(\Lambda^k V)$ by the rule 
$$P(U)={f_1\wedge \dotsb \wedge f_k}.$$
The point $P(U)$ has the following form in the basis $\{e_{i_1}\wedge \dotsb \wedge e_{i_k} \}_{i_1< \dotsb <i_k} $ of $\La^k V$:
$$P(U)= {\sum_{i_1< \dotsb <i_k}p_{i_1 \dotsc i_k} e_{i_1}\wedge \dotsb \wedge e_{i_k}}.$$
Then the homogeneous coordinates $p_{i_1 \dotsc i_k}$ of  $P(U)$ are called the Plucker coordinates of $U$,
$P(U)$ is uniquely determined by   $U$, and the following assertions hold.

\begin{enumerate} 
 \item[{\sffamily(i)}]    
The subset of all points  $p\in \mathbb{P}$$(\Lambda^k V)$ of the form $p=P(U)$   is  closed in  $\mathbb{P}$$(\Lambda^k V)$; this subset $G(k,n)$ (the  Grassmannian or Grassmann variety) is defined by the relations   
\begin{equation} \label{E:d}
      \sum_{r=1}^{k+1}(-1)^r p_{i_1 \dotsc i_{k-1} j_r}p_{j_1 \dotsc \widehat{j_r} \dotsc j_{k+1}}=0
\end{equation} 
for all sequences   $i_1\dotsc i_{k-1}$ and $j_1 \dotsc j_{k+1}$.

\item[{\sffamily(ii)}]
  $\dim G(k,n)=k(n-k)$.
\item [{\sffamily(iii)}]  $G(k,n)$ is irreducible {\normalfont(see  Section~2.2.7 of~\cite{VO95}  )}.

\item[{\sffamily(iv)}] Suppose, for example, that  $p_{1\dotsc k} \ne 0$. If  $p=(p_{i_1\dotsc i_r})=P(U)$, then  $U$ has a basis $\{ \itd f k \}$ such that   
\begin{equation} \label{E:e}
 f_i=e_i+\sum_{r>k} a_{ir}e_r \qquad \mbox{for }  i\odo k,
\end{equation}
where
\begin{equation} \label{E:f}
  a_{ir}=(-1)^{k-i} \frac{p_{1\dotsc \widehat{i}\dotsc kr}}{p_{1\dotsc k}}.
\end{equation}
(Here $\widehat i $ means that the index $i$ is discarded.)
\end{enumerate}

\end{utv}

\begin{utv} \label{T:v2} 
{\normalfont (see~\cite[Seciton 14.7]{F89} or~\cite[Section 1.5]{G82} ) }
Concider an $n$-dimensional vector space $V$. For any increasing sequence $$
0 \subsetneq V_1\subsetneq \dotsc \subsetneq V_k 
$$
of subspaces of $V$
 put
\begin{equation*}
W(\itd V k)=\{p=P(U)\in G(k,n) \mid \dim (U\cap V_i)\geq i  \mbox{ for } i\odo k \}.
\end{equation*}
Then all the subsets $W(\itd V k)$, called Schubert cells, are closed in $G(k,V)$, and  
\begin{equation} \label{E:Y2}
\dim W(\itd V k)= \sum_{i=1}^k (  \dim V_i-i). 
\end{equation}
\end{utv}

\begin{sled}
Let $V$ be a vector space of dimension $n$, and let ${s_0=\max \{ 0, k+m-n\}}$. Given an $m$-dimensional subspace ${U\subset V}$ $(m>0)$, put
\begin{equation} \label{E:Z2}
 G_s=G_s(U,V,k)=\{P(U')\in G(k,n) \mid \dim U \cap U' \geq s  \} 
\end{equation}
 for $ s=s_0, s_0+1,\dotsc, \min \{k,m\}$.
Then the sets $G_s$ are closed in $G(k,n)$, and
 \begin{equation} \label{E:L2}
\dim G_s = s(m-s)+ (k-s)(n-k). 
\end{equation}
\end{sled}
\doc
It is easy to check that the sets $G_s$ are Schubert cells. Indeed, choose a basis $\{\itd e n\}$ of $V$ such that $\itd e m$ span $U$. Then    ${G_s=W(\itd V k)}$ where
\begin{alignat*}{2}
&     V_i= \langle \itd e {m-s+i}\rangle\ \ \    &&\mbox { for } i \odo s, \\
&V_i = \langle \itd e {n-k+i}\rangle   &&\mbox{ for } i=s+1,\dotsc, k.
\end{alignat*}
(The conditions $\dim(U'\cap V_i)\geq i$ for $i> s$ are trivial , and for $i\leq s$ they follows from the condition
 $\dim (U' \cap V_s)=\dim(U'\cap U) \geq s.$)  Using~(\ref{E:Y2}), we get
$$
\dim G_s= \sum_{i=1}^s(m-s)+\sum_{i=s+1}^k(n-k)=s(m-s)+(k-s)(n-k) .\ \bl
$$

\section{Proofs} \label{S3}
\begin{lem} \label{tj}
Let $K$ be an infinite field. Denote by  $\bar K$ its algebraic closure.
Suppose that  a generic Lie algebra of  ${\mathfrak N}_2^{\bar K}(n,t)$ has property $\ma$ and for any Lie algebra  $L\in \nilk$ without property $\ma$ the Lie algebra $\bar L=L_K\otimes \bar K \in {\mathfrak N}_2^{\bar K}(n,t)$ does not have property $\ma$ too. Then $\ma$ is also a generic property of $\nilk$.
\end{lem} 
\doc 
The $t$-tuples  $\Ph(L)$ and $\Ph(\bar L)$, corresponding to the Lie algebras 
$L\in \nilk$   and $\bar L  \in {\mathfrak N}_2^{\bar K}(n,t)$ respectively, can be defined by the same $t$-tuple of matrices with coefficients in $K$.
So, under the  condition of the Lemma, the set of point of $M_{n,t}(K)$ corresponding to the algebras of $\nilk$ without property $\ma$ satisfies some finite system of homogeneous equations, say $(\ast)$, over $\bar K$. If the variables take values in $K$, $(\ast )$  is equivalent to some finite system of homogeneous  equations over $K$. The last system is not zero on $M_{n,t}(K)$ identically, since $K$ is infinite and $(\ast)$ defines proper subset of $M_{n,t}(\bar K)$.~$\bl$\\

It follows from Lemma~\ref{tj} that it is enough to prove Theorems~\ref{ti}-\ref{tf} just for the case of algebraically closed ground field $K$. So,  from now on  if not stated otherwise  we will assume that  $K$ is algebraically closed.

Propositions~\ref{ta} and~\ref{tb} show that we can define the correspondence between the elements of $\nil$ and $G(t,B_n )$ by the rule
 $\ph(L)=\ph_1\wedge \dots \wedge \ph_t$ for $L\in\nil$ if $\Ph(L)=(\itd \ph t)$, and vice versa,  $L(\ph)=L(\itd \ph t)$ for $\ph=\ph_1\wedge \dots \wedge \ph_t\in G(t, B_n)$. 
 
\begin{lem} \label{tc}
$\ma$ is a generic property of $\nil$ if and only if the set of  points of $G(t, B_n)$ corresponding to the algebras without property $\ma$ belongs to some proper Zariski-closed subset of $G(t, B_n)$.
\end{lem} 

\doc 
Consider a regular map $f:M_{n,t}\to G(t, B_n)$ such that $f$ takes each $t$-tuple of alternating bilinear forms to the  subspace spanned by it, that is, 
$$
f(\itd \ph t)= \ph_1\wedge \dots \wedge \ph_t.
$$

By definition, Lie algebras $L(\Ph)$ and $L(f(\Ph))$ corresponding to the points $\Ph \in M_{n,t}$ and $f(\Ph)$ respectively are isomorphic.

Denote by $M(\ma)\subset M_{n,t}$ and $G(\ma)\subset G(t,B_n)$ the sets of points corresponding to the algebras without property $\ma$. Clearly, $f^{-1}(G(\ma))=M(\ma)$. 
Suppose that $G(\ma)$  belongs to some proper closed subset $X\subsetneq G(t,B_n)$. Then $M(\ma)$ is contained in  $f^{-1}(X)$.  Since $f$ is surjective and continuous, $f^{-1}(X)$ is a proper closed subset of $M_{n,t}$.

Inversely, suppose that $M(\ma)$ belongs to some proper closed subset $Y \subsetneq M_{n,t}$.
A fiber $f^{-1} (\ph)$ consists of all the bases of a given $t$-dimensional vector space considered up to proportionality. It can be parametrized by invertible  matrices of size $t\times t$. Hence, all the fibers $f^{-1} (\ph)$ for $\ph \in G(t,B_n)$ are isomorphic and of the same dimension $t^2-1$.

Obviously, $L(\Ph)$ does not have property $\ma$ only if $f^{-1} (f(\Ph))\subset Y$. Hence, we can consider a set $Y_0= \{\Ph \in Y \mid  \dim_\Ph (f^{-1} (f(\Ph))\cap Y)= t^2 -1\}$ instead of $Y$. An application of  Proposition~{\ref{tl}} to the map  $\bar f = f|_Y$ yelds that $Y_0$ is also closed.

Denote by $X_0$ the closure of $f(Y_0)$ in  $G(t,B_n)$. We have $G(\ma)\subset X_0$. Let us show that $X_0$ is proper subset of  $G(t,B_n)$. Indeed, by~(\ref{ez}), for any irreducible component $Y'_0 \subset Y_0$ we have
$$
\dim Y'_0=\dim X'_0 +t^2 -1,
$$
where  $X'_0 \subset X_0$  is a closure of the image $f(Y'_0)$.
Recalling that $\dim M_{n,t} = \dim G(t,B_n) +t^2-1$ and $\dim Y'_0 < \dim M_{n,t}$, we obtain  $\dim X'_0 < \dim G(t, B_n)$.
Consequently,  $\dim X_0 < \dim G(t,B_n)$.~$\bl$

\subsection{Proof of Theorem~\ref{ti}} 

Consider a Lie algebra $L \in \nil$. Let $S=[L,L]$, $L=V \oplus S$    and $\Ph(L)=(\itd \ph t)$. Fix  a basis $\itd e n$ for  $V$ and identify forms $\ph_i$ with    their matrices with respect to this basis.
$Z(L)$ is strictly greater than  $S$ if and only if there exists nonzero element $c \in V$  such that
\begin{equation} \label{es}
\ph_i(c, e_j)=0 \qquad   i\odo t,\ j\odo n.
\end{equation}

 Denote by $D$ the subset of $H={\mathbb P}(V) \times {\mathbb P}(B_n^t)$ consisting of pairs $(c,\ph)$ satisfying~(\ref{es}). The system of equations~(\ref{es}) is linear in each set of variables $c$ and $\ph$. Therefore, by Proposition~\ref{T:e2}, $D$ is a projective variety.

Consider the projection $\pi: D \to {\mathbb P}(B_n^t)$. Since $ \pi$ is a regular map, $ \pi(D)$ is closed in $ {\mathbb P}(B_n^t)$. Moreover,  $\pi(D)\cap M_{n,t}$ consists of $t$-tuples  $(\ph_1,\dotsc,\ph_t)$ corresponding to the algebras with non-trivial center modulo the commutator ideal.
Obviously,  $\pi(D)\cap M_{n,t}$ is a proper subset except the case when $t=1$ and $n$ is odd  in which the only form $\ph_1$ defining a Lie algebra is always degenerate.  In the last case the dimension of the center of a Lie algebra is greater than 2 only if the rank of  $\ph_1$ is strictly less than $n-1$. This condition defines a  proper closed subset in $M_{n,1}$.~$\bl$

\subsection{Proof of Theorem~\ref{ts}}\label{s6}
\emph{Proof of Theorem~\ref{ts}. }
The Main Lemma of~\cite{M7} states that if the positive integers $t\geq 2 $, $k$ and   $n$ satisfy the inequality 

\begin{equation} \label{et}
2n \geq t(k-1)+2k,   
\end{equation} 
then  for any $t$-tuple $\{\itd \ph t\}\in B_n^t(K)$ there exists a $k$-dimensional subspace that is simultaneously isotropic for all of the forms $\itd \ph t$, i.e. on which all the forms are zero.

It follows  easily from the proof of the Main Lemma that the set $\pi_2(S)$ (in the notation of~{Section 2.2} of~\cite{M7}) consisting of all the $t$-tuples of forms with common $k$-dimensional isotropic subspace is a Zarisski-closed subset of $\Pm(B_n^t)$. Furthermore, we have $\pi_2(S)=\Pm(B_n^t)$ if  $t>1$ and the inequality~(\ref{et}) holds, and  $\pi_2(S)$ is a proper subset if~(\ref{et}) is not  true.  

Let $L(\Ph)=V\oplus S\in \nil$  be a Lie algebra associated with a tupl $\Ph$.  A subalgebra $H\leq L$ is commutative if and only if the subspace $H/(H\cap S)$ is isotropic for all of the forms of $\Ph$. Therefore each algebra of $\nil$ contains a commutative subalgebra of dimension $s=k+t$, where $k$ is the maximal integer satisfying~(\ref{et}), that is $k=\left[\frac {2n+t}{t+2}\right]$. And a generic Lie algebra of $\nil$ has no abelian subalgebras of dimension greater than $s$.~$\bl$\\

If the ground field $K$ is not algebraically closed, then  
the maximal dimension of commutative subalgebras of a Lie algebra $L$ over $K$ may be strictly less than that of $\bar {L}=L_K\otimes \bar  K$ over $\bar K$.
And so the condition on the ground field to be  closed in the first part of Theorem~\ref{ts} is essential.

\begin{exam}
Let $K$ be a subfield of the  field of real numbers. There exists a Lie algebra $L\in \mathfrak{N}^K_2(4,3)$ without commutative subalgebras of dimension 5.
\end{exam}
\doc Concider a Lie algebra ห่ $L=V\oplus S$ defined by the following 3-tuple of matrices
\begin{equation*}
\ph_1=
\begin{pmatrix}
0&1&0&0\\
-1&0&0&0\\
0&0&0&1\\
0&0&-1&0
\end{pmatrix}
,\ \ph_2=
\begin{pmatrix}
0&0&1&0\\
0&0&0&-1\\
-1&0&0&0\\
0&1&0&0
\end{pmatrix}
,\ \ph_3=
\begin{pmatrix}
0&0&0&1\\
0&0&1&0\\
0&-1&0&0\\
-1&0&0&0
\end{pmatrix}
.
\end{equation*}
Elements $x,y\in V $ commute if and only if $\ph_i(x,y)=0$ for  $i=1,2, 3.$ Hence the set of all elements $y\in V$ commuting with a given element $x$ with coordinates $(x_1,x_2,x_3,x_4)$
is a solution space of the system of homogeneous linear equations with  the matrix
\[
  M(x)= 
\begin{pmatrix}
x_2&-x_1&x_4&-x_3\\
x_3&-x_4&-x_1&x_2\\
x_4 &x_3& -x_2&-x_1 
\end{pmatrix}.
\]
The system has a solution  not proportional to x if and only if $M(x)$ is rank deficient.
Let us compute $3\times 3$ minors of $M(x)$:
\begin{equation*}
 A_i=(-1)^i x_i(x_1^2+x_2^2+x_3^2+x_4^2) \qquad \mbox{for }i=1,2,3,4.
\end{equation*} 
Obviously, if all   $A_i=0$ and  $x_i \in \mathbb{R}$, then $x_1=x_2=x_3=x_4=0$. Consequently, $L$ has no commutative subalgebras of dimension greater than 1 modulo commutator ideal.~$\bl$

\subsection{Proof of Theorem~\ref{tf}}
For any subspace $U$ of an $n$-dimensional vector space $V$, put
\begin{equation*} \label{eon}
N_0(U)=\{\ph \in B_n \mid \ph(x,u)=0 \mbox{ for } x\in V \mbox{ and } u\in U\}.
\end{equation*}
Choose a basis $\itd e n $ of $V$ such that $\itd e k$ span $U$. Now an alternating bilinear form $\psi$ lies in $N_0(U)$ if and only if the matrix of $\psi$ is zero except  lower right-hand $(n-k) \times (n-k)$ submatrix  in this basis. Thus,
\begin{equation} \label{eh}
\dim N_0(U)=\frac {(n-k)(n-k-1)}2.
\end{equation}     

\begin{lem} \label{te}
Consider a Lie algebra $L \in \nil$. Let $S=[L,L]$, $L=V \oplus S$    and $\Ph(L)=(\itd \ph t)$.  $L$  has a property $\ms$$(n_0,t_0)$
if and only if there exists a subspace $U \subset V$ of dimension $n-n_0$ such that
\begin{equation} \label{ec}
 \dim (N_0(U) \cap  \langle \itd \ph t \rangle) \geq t_0.
\end{equation}
\end{lem}
\doc L has property $\ms$$(n_0,t_0)$ if and only if there exists an ideal $I\triangleleft L$ of dimension $(n+t-n_0-t_0)$ such that $\dim (I\cap S)= t-t_0$. Indeed, in this case we have $L/ I \in {\mathfrak N}_2(n_0,t_0)$, since the image of a commutator ideal is the commutator ideal of the image. And so we can take $I$ as a kernel of a desired surjective map.

Assume that such $I$ exists. Choose a basis $\itd {z'} t$  of $S$ such that $ z'_{t_0+1}, \dots, z'_t\in (S\cap I)$. Let $(\itd {\ph'} t)$ be a $t$-tuple of forms associated with $L$ in this basis. Denote by $U$ the projection of $I$ on $V$ along $S$. Clearly, $\dim U =n-n_0$. Since for any $x\in V$ and $u\in U $ we have $[x,u]=\sum_{i=1}^t \ph_i'(x,u)z'_i \in \langle z_{t_0+1}, \dots, z_t \rangle$, we get $\itd {\ph'}{t_0}\in N_0(U)$. It follows from~(\ref{eb}) that $\itd {\ph'}{t_0}$ are linear combinations of $\itd \ph t$ and, by Proposition~\ref{ta}, they are linearly independent. Consequently,~(\ref{ec}) holds.

On the contrary,     let $U$ be as in the lemma. Choose $t_0$ linearly independent forms 
 $\itd  {\ph'}{t_0}\in (N_0(U) \cap  \langle \itd \ph t \rangle)$ and extend  them to a basis  $\itd {\ph'} t$ of  $\langle \itd \ph t \rangle$. It follows from~(\ref{ea}) and~(\ref{eb}) that we can choose a basis $\itd {z'} t$ of $S$ such that $\Ph(L)=(\itd  {\ph'} t)$ in this basis. For any  $x\in V$ and  $u\in U $ we have $[x,u]=\sum_{i=t_0+1}^t \ph'_i(x,u)z'_i \in \langle z_{t_0+1}, \dots, z_t \rangle$. Hence, the linear span $I=\langle U, z_{t_0+1}, \dots, z_t\rangle$ is a desired ideal.~$\bl$

\begin{lem} \label{tm} 
Suppose that the positive integers $n,\ n_0,\ t$ and $ t_0$ satisfy the inequality~(\ref{eq}) and let  $n\geq n_0 $. Then any Lie algebra $L\in \nil$ has the property $\ms$$(n_0,t_0)$.
\end{lem}
\doc
In the notation of Lemma~\ref{te}, for any Lie aglebra $L\in \nil$ and for any subspace $U\subset V $ of dimension $n- n_0$ we have
\begin{align*}
\dim (N_0(U) \cap \langle \itd \ph t \rangle  )&\geq \dim N_0(U) +\dim\langle \itd \ph t \rangle -\dim B_n = \\
 &\overset{(\ref{eh})}{=}
\frac{n_0(n_0-1)}2+ t - \frac{n(n-1)}2.
\end{align*}
Combining this with inequality~(\ref{eq}), we get~(\ref{ec}).~$\bl$\\

Obviously, if  $n < n_0 $ or $t<t_0$, there are no Lie algebras in $\nil$ having the property $\ms$$(n_0,t_0)$. And so in what follows to prove the theorem we fix integers $n,\  n_0,\  t,\ t_0,\  k$ such that $n=n_0+k,$ $t\geq t_0\geq 1,$ $  n_0\geq 2,$ $ k\geq 0$ and  an $n$-dimensional vector space $V$.
Denote by $D$ the following  subset of the direct product 
 $ H=G(k,V)\times G(t,B_n) \subset {\mathbb P}(\Lambda^k    {V}) \times {\mathbb P}(\La^t B_n)$
$$
D=\{(p,\ph)\in H \mid  p=P(U),\ \ph=P(\Omega),\  \dim (N_0(U) \cap  \Omega) \geq t_0 \}.
$$

\begin{lem} \label{T:p2}
$D$ is a projective variety.
\end{lem}
\doc
 Fix a basis $\itd e n $ of $V$ and a basis $\itd  \psi {\frac{n(n-1)}2}$ of $B_n$. 
Then the coordinates $(p_{i_1 \dotsc i_k})$ and $(\ph_{j_1 \dots j_t})$ on ${\mathbb P}(\Lambda^k    {V}) $ and $ {\mathbb P}(\La^t B_n)$  respectively (and, hence, on $H$) are naturally defined.

Let us prove that $D$ is closed in $H$. For this we  consider the covering of the projective variety $H$ by the open sets $${O_{i_1 \dotsc i_k, j_1\dots j_t}=\{ (p,\ph)\in H | p_{i_1 \dotsc i_k} \ne 0 \mbox{ and } \ph_{j_1 \dots j_t} \ne 0  \}}$$ and verify that $O_{i_1 \dotsc i_k, j_1\dots j_t} \cap D$ is closed in~$O_{i_1 \dotsc i_k, j_1\dots j_t}$ for any sequences $i_1 \dotsc i_k$ and  $j_1\dots j_t$ with $i_1 <\dotsb <i_k,$ and $ j_1<\dotsb <j_t$. Then, by Proposition~\ref{T:a2}, $D$ is closed in $H$. Hence, $D$ is also a projective variety.

Let us consider an arbitrary point $(p,\ph)\in G(k,V)\times G(t,B_n)$. We can assume without loss of generality that $p_{1\dotsc k} \ne 0$ and $\ph_{1\dots t }\ne 0$. Let 
${p=P(U)}$ and  $\ph=P(\Omega)$. Then, according formula~(\ref{E:e}), $U$ has the basis

\begin{equation} \label{eg}
 f_i=e_i+\sum_{r=k+1}^n a_{ir}e_r \qquad \mbox{for }i\odo k,
\end{equation}
where
\begin{equation} \label{ed}
  a_{ir}=(-1)^{k-i} \frac{p_{1\dotsc \widehat{i}\dotsc kr}}{p_{1\dotsc k}};
\end{equation}
and  $\Omega$  has the basis 
\begin{equation} \label{ef}
 \ph_i=\psi_i+\sum_{r=t+1}^{\frac{n(n-1)}2} b_{ir}\psi_r \qquad \mbox{for }i\odo t,
\end{equation}
where
\begin{equation} \label{ee}
  b_{ir}=(-1)^{t-i} \frac{\ph_{1\dotsc \widehat{i}\dotsc tr}}{\ph_{1\dotsc t}}.
\end{equation}

It follows from~(\ref{eg}) that the vectors $\itd f k, e_{k+1}, \dots, e_n$  form a basis of $V$.
Denote by $C$ the transformation matrix from the basis $\itd e n$ to this new one and by  $\itd A {\frac{n(n-1)}2}$ the matrices of  $\itd \psi {\frac{n(n-1)}2}$ in the initial basis. Then, by formula~(\ref{ef}), matrices of $\itd \ph t$ have the form
$$
F_i=C^\bot (A_i  +   \sum_{r=t+1}^{\frac{n(n-1)}2} b_{ir}A_r) C  \qquad \mbox{for }i\odo t,
$$
in the new basis and, as we mentioned above, $N_0(U)$ consists of the forms whose matrices are zero except  lower right-hand $(n-k) \times (n-k)$ submatrix. 

Let $\itd E {\frac {n_0(n_0-1)}2}$ be a basis of $N_0(U)$. (These are matrices with constant coefficients.) The condition $\dim (N_0(U) \cap  \Omega) \geq t_0$  holds if and only if the rank of the  vector system $\{\itd F t, \itd E {\frac {n_0(n_0-1)}2}\}$ is not greater than $s= \frac {n_0(n_0-1)}2+t-t_0$. The last condition in turn is equivalent to the fact that all the minors of size $s+1$ in the corresponding matrix are zero. That is, we  have a system of polynomial equations in variables $a_{ir}$ and $b_{ir}$. Replacing this variables using~(\ref{ed}) and~(\ref{ee}) and multiplying both sides of the equations by the appropriate powers of  $p_{1\dots k}$  and $\ph_{1\dots t}$, we obtain the system of equations homogeneous separately in each set of variables $p$ and $\ph$. According to Proposition~\ref{T:e2}, the set $D \cap O_{1\dotsc k, 1\dots t}$, defined by the last system, is closed in $O_{1\dotsc k, 1 \dots t}$.~$\bl$\\

Now consider the projections 
 $\pi_1:D \to {\mathbb P}(\Lambda^k   {V})$ and $\pi_2: D \to {\mathbb P}(\Lambda^t B_n)$ such that ${\pi_1(p,\ph)=p}$, ${\pi_2(p,\ph)=\ph}$. These are regular maps. Clearly  $\pi_1(D)=G(k,V)$.

\begin{lem} \label{tu}
$\ph_1\wedge \dotsc\wedge\ph_t \in G(t, B_n)$ belongs to $\pi_2(D)$ if and only if the corresponding Lie algebra $L(\itd \ph t)$ has the property $\ms$$(n_0,t_0)$. $\pi_2(D)$ is closed in  $G(t, B_n)$.
\end{lem}
\doc The lemma follows from Lemma~\ref{te} and from the fact that the image of the projectiv variety $D$  under the regular map $\pi_2$ is closed.~$\bl$ \\

\begin{lem} \label{tw}
If the integers $n,\ n_0,\ t$ and $ t_0$ do not satisfy the relation~(\ref{eq}), then for any point $p \in G(k,V)$  the fiber $\pi_1^{-1}(p)$ is an irreducible projective variety  of dimension
\begin{equation} \label{el}
\dim \pi_1^{-1}(p) =  t_0\left(\frac{n_0(n_0-1)}2- t_0\right)+ (t-t_0)\left(\frac {n(n-1)} 2-t \right).
\end{equation}
\end{lem}
\doc For any point $p=P(U)\in G (k,V)$ we have
$$
\pi_1^{-1}(p)=\{P(\Omega)\in G(t,B_n) \mid \dim (N_0(U) \cap \Omega) \geq t_0  \}.
$$
That is, the fiber $\pi_1^{-1}(p)$ is a set of type~(\ref{E:Z2}):
$$
\pi_1^{-1}(p)=G_{t_0}(N_0(U),B_n,t).
$$
Since~(\ref{eq}) doesn't holds, $t_0$ satisfies the condition of the Corollary of Proposition~\ref{T:v2}. We apply this corollary,  taking in account that, by~(\ref{eh}), $\dim N_0(U)=\frac {n_0(n_0-1)}2$ and  $\dim B_n=\frac {n(n-1)}2 $, and thus  obtain~(\ref{el}).~$\bl$

\begin{lem} \label{tx}
$D$ is an irreducible variety of dimension
\begin{align} \label{ej}
\dim D = 
 t_0\left(\frac{n_0(n_0-1)}2- t_0\right)+ (t-t_0)\left(\frac {n(n-1)} 2 -t \right) + n_0(n-n_0).
 \end{align}
\end{lem}
\doc
By the previous lemma, all the fibers $ \pi_1^{-1}(p)$ are irreducible and of the same dimension.
The image $\pi_1(D)=G(k,V)$ is also irreducible by Proposition~\ref{T:u2}~(iii). So, by Proposition~\ref{T:A2}, $D$ is irreducible and we can use formula~(\ref{ez}) to compute its dimension. Namely, 
$\dim D = \dim \pi_1^{-1}(p) + \dim G(k,V)$. Using Lemma~\ref{tw} and Proposition~\ref{T:u2}~(ii), we get~(\ref{ej}).~$\bl$\\

\emph{Proof of Theorem~\ref{tf}.}  
To prove the first part of the theorem it is enough to show that under the condition~(\ref{er}) the indequality  $\dim \pi_2(D) < \dim G(t, B_n)$ holds and apply Lemmas~\ref{tc} and~\ref{tu}. By~(\ref{ez}), this inequality follows from 
$$\dim D < \dim G(t, B_n).$$
  Substituting for $\dim D$ using Lemma~\ref{tx} and recalling that $\dim G(t, B_n)= t\left( \frac{n(n-1)}2 -t \right)$ we obtain the relation equivalent to~(\ref{er}).

The second part of the theorem follows from Lemma~\ref{tm}. So, Theorem~\ref{tf} is proved.~$\bl$\\

\emph{Proof of the Corollary of Theorem~\ref{tf}.}
A Lie algebra $L$ does not admit a surjective homomorphism on a non-abelian Lie algebra of dimension $N$, if it does not have the property $\ms (n_0, t_0)$ for any positive integers $n_0, t_0 $ such that $n_0+t_0=N$ and $t_0\leq \frac{n_0(n_0-1)}2$. By Theorem~\ref{tf}, this condition holds for almost all the Lie algebras $L \in \nil$ if for all of the stated values of $n_0$ and $t_0$ we have~(\ref{er}), that is,

\begin{align*}
t<\frac{n^2}2+ \frac 1 2 n -2n_0 - \frac{N(n-N)}{N-n_0}- \frac{n_0(n_0-1)}2. \label{em}
\end{align*}
If $N\leq n$, the right term of the inequality is decreasing in $n_0$ function and, hence, it takes the minimum value at  $n_0=N-1$.~$\bl$
 
\textbf{Acknowledgment.} The author would like to express her gratitude to Professor A.Yu.~Olshanskii  for   suggestion of the problem.

\small \emph{E-mail address:} mariamil@yandex.ru
\end{document}